\documentclass[12pt]{article}

\usepackage{amsmath}
\usepackage{amsthm}
\usepackage[pdftex]{graphicx}
\usepackage{ccaption}

\newtheorem{theorem}{Theorem}[section]

\newtheorem{conjecture}[theorem]{Conjecture}

\newcommand{\script}[1]{\text{$\cal{#1}$}}
\newcommand{\dist}[0]{dist}
\newcommand{\greg}[1]{}

\theoremstyle{remark}
\newtheorem*{remark}{Remark}

\theoremstyle{definition}
\newtheorem{definition}[theorem]{Definition}

\usepackage[margin=1in]{geometry}        

  \captionnamefont{\bfseries}
  \captiontitlefont{\small\sffamily}
  \captiondelim{ --- }
  \hangcaption

\begin{document}

\title{Contributions to Seymour's Second Neighborhood Conjecture}
\author{James N. Brantner, Greg Brockman,\\ Bill Kay, and Emma E. Snively\\
\small Erskine College, Harvard University,\\
\small University of South Carolina, Rose-Hulman Institute of Technology
}

\maketitle
\begin{abstract}
Let $D$ be a simple digraph without loops or digons (i.e. if $(u,v)\in E(D)$, then $(v,u)\not\in E(D)$). For any $ v\in V(D)$ let $N_1(v)$ be the set of all vertices at out-distance 1 from $v$ and let $N_2(v)$ be the set of all vertices at out-distance 2.  We provide sufficient conditions under which there must exist some $v \in V(D)$ such that $|N_1(v)| \leq |N_2(v)|$, as well as examine properties of a minimal graph which does not have such a vertex. We show that if one such graph exists, then there exist infinitely many strongly-connected graphs having no such vertex. 
\end {abstract}

\section{Introduction}

For the purposes of this article, we consider only simple nonempty digraphs (those containing no loops or multiple edges and having a nonempty vertex set), unless stated otherwise.  We also require that our digraphs contain no digons, that is, if $D$ is a digraph then $(u,v) \in E(D) \Rightarrow (v,u)\notin E(D)$.  If $i$ is a positive integer, we denote the $i^{\text{th}}$ neighborhood of a vertex $u$ in $D$ by $N_{i,D}(u) = \{v \in V(D) | \dist_D(u, v) = i\}$, where $\dist_D(u, v)$ is the length of the shortest directed path from $u$ to $v$ in $D$ (if there is no directed path from $u$ to $v$, we set $\dist_D(u, v)=\infty$).  If $D$ is clear from context, we simply write $N_i(u)$ and $dist (u,v)$. We also may wish to consider the $i^\text{th}$ in-neighborhood of a vertex $N_{-i}(u) = \{v\in V(D) | \dist(v, u) = i\}$. In addition, if $V' \subseteq V(D)$, we let $D[V']$ be the subgraph of $D$ induced by $V'$. 
    
Graph theorists will be familiar with the following conjecture due to Seymour (see~\cite{dean}), now more than a decade old:

\begin {conjecture}[Seymour's Second Neighborhood Conjecture]
\label{SeymourConjecture}
Let $D$ be a directed graph. Then there exists a vertex $v_0 \in V(D)$ such that $|N_1(v_0)|\leq |N_2(v_0)|$.  
\end{conjecture}

In 1995, Dean~\cite{dean} conjectured this to be true when $D$ is a tournament.  Dean's Conjecture was subsequently proven by Fisher~\cite{fisher} in 1996.  Further, in their 2001 paper Kaneko and Locke~\cite{kaneko} showed Conjecture \ref{SeymourConjecture} to be true if the minimum outdegree of vertices in $D$ is less than 7, and Cohn, Wright, and Godbole~\cite{godbole} showed that it holds for random graphs almost always.  And finally, in 2007 Fidler and Yuster~\cite{yuster} proved that Conjecture \ref{SeymourConjecture} holds for graphs with minimum out-degree $|V(D)| - 2$, tournaments minus a star, and tournaments minus a sub-tournament. While over the years there have been several attempts at a proof of Conjecture \ref{SeymourConjecture}, none of these have yet been successful.

For completeness, we introduce the related Caccetta-H\" aggkvist conjecture~\cite{caccettahaggvist}, which was posed in 1978:

\begin{conjecture}[Caccetta-H\" aggkvist Conjecture]
If $D$ is a directed graph with minimum outdegree at least $|V(D)|/k$, then $D$ has a directed cycle of length at most $k$.
\label{CH}
\end{conjecture}

Conjecture \ref{SeymourConjecture} would imply the $k=3$ case of Conjecture \ref{CH}.  Much work has been done on Conjecture \ref{CH}, including an entire workshop in 2006 sponsored by AIM and the NSF, yet Conjectures \ref{SeymourConjecture} and \ref{CH} both remain open.

We do not seek to prove Conjecture \ref{SeymourConjecture} in this paper.  Rather, we prove the conjecture for various classes of graphs.  We then take a different tack and provide conditions that must be satisfied by any appropriately-defined minimal counterexample to Seymour's Second Neighborhood Conjecture. This provides tools with which the conjecture can be approached; in one direction it may aid in showing the nonexistence of such a graph, while in the other direction we restrict the search space of possible counterexamples.

\section{Definitions}

We begin our investigation by defining some useful terms.

\begin{definition} 
\label{Seymour}
Suppose that $D$ is digraph and $u\in V(D)$.  We say that $u$ is \textit{satisfactory} if $|N_1(u)| \leq |N_2(u)|$.  Also, $u$ is a \textit{sink} if $|N_1(u)|=0$.  Note that a sink is trivially satisfactory.
\end {definition}

\begin{definition}
\label{minimal criminal}
Let $\script{A} = \{D | \textit{$D$ is a simple directed graph with no satisfactory vertices}\}$ be the set of counterexamples to Seymour's Second Neighborhood Conjecture.  Let $\script{A'} = \{D\ |\ |E(D)| = \min_{H\in \script{A}} |E(H)| \}$ be the set of graphs in $\script{A}$ with the fewest number of edges.  Finally, let $\script{A''} = \{D\ |\ |V(D)| = \min_{H\in \script{A'}} |V(H)| \}$ be the set of graphs in $\script{A'}$ with the fewest number of vertices.  We will refer to any element of $\script{A''}$ as a \textit{minimal criminal}.  Note that $\script{A''}$ is empty if and only if Conjecture \ref{SeymourConjecture} is true.
\end{definition}

\begin{definition}
\label{walkable neighborhood}
Let $D$ be a digraph.  Suppose that $u \in V(D)$.  We define $W_D(u)=\{v | \dist(u,v) \neq \infty\}$ to be the \textit{walkable neighborhood} of $u$ with respect to $D$. If $D$ is clear from context, we simply write $W(u)$.

Also define $A_{s,D}(u)= |N_1(u)| - |N_2(u)|$ to be the \textit{anti-satisfaction} of $u$.  As usual, if $G$ is clear from context, we simply write $A_s(v)$.  Notice that $u$ is satisfactory if and only if $A_s(u) \leq 0$.
\end {definition}

\begin{definition}
Again let $D$ be a directed graph.  Recall that a transitive triangle $T$ is a directed graph on three nodes $a,b,c$ such that $(a,b), (a,c), (b,c)\in E(T)$.  If $(u,v) \in E(D)$, we say that edge $(u,v)$ is the \textit{base of a transitive triangle} if $u$ and $v$ share a common first neighbor; that is, $|N_1(u)\cap N_1(v)| \geq 1$.

\begin{figure}[htb]
\centering
\includegraphics{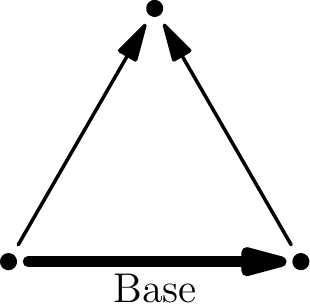}
\caption{Demonstration of an edge that is the base of a transitive triangle}
\label{fig:transitive triangle}
\end{figure}

If, for distinct $t,u,v,w \in V(D)$, we have that $(t,u),(u,w),(t,v),(v,w)\in E$ then we call $\{(t,u),(u,w),(t,v),(v,w)\}$ a \textit{2-directed diamond}. We say the edges $(t,u), (t,v)$ are the \textit{bases} of the 2-directed diamond.

\begin{figure}[htb]
\centering
\includegraphics{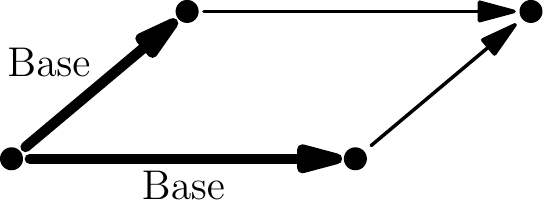}
\caption{Demonstration of the bases of a 2-directed diamond}
\label{fig:seymour diamond}
\end{figure}
\end{definition}

We now have the tools to delve into our results.

\section{Directed cycles and underlying girth}

In this section we show that certain classes of graphs satisfy Seymour's Second Neighborhood Conjecture.  The following theorem shows that directed cycles are necessary for a graph to be a counterexample to the conjecture.

\begin{theorem}
\label{D.Cycles}
If a digraph contains no directed cycles, then it must have a satisfactory vertex.
\end{theorem}
\begin{proof}
Let $D$ be a directed graph, and suppose that $D$ contains no satisfactory vertices.  Then $D$ has no sink, as noted in Definition \ref{Seymour}.  It is a well-known fact that a graph with no sinks has a directed cycle.  We include the standard proof, however, since the same technique will be useful to us later: Because $D$ has sink, for $v\in V(D)$, $|N_1(v)| > 0$.  Pick an arbitrary vertex $v_0 \in V(D)$, and consider the sequence $\{v_i\}_{i=0}^{|V(D)|}$ defined recursively by $v_{i+1}\in N_1(v_{i})$ for $i\geq 0$.  By the Pigeonhole principle, there exist some $r\neq s$ such that $v_r = v_s$.  Then we note that the sequence of edges $(v_r, v_{r+1}), (v_{r+1}, v_{r+2}), \ldots, (v_{s-1}, v_s = v_r)$ defines a dicycle in $D$, thus completing our proof.
\end{proof}

The following theorem provides another sufficient condition for a graph to contain a satisfactory vertex:

\begin{theorem}
\label{girth}
Let $D$ be a directed graph containing no transitive triangles.   Then $D$ contains a satisfactory vertex.
\end {theorem}

\begin{proof}
Let $v_0\in V(D)$ have the minimal out-degree in $D$.  If $|N_1(v_0)| = 0$, then $v_0$ is a sink and hence a satisfactory vertex.  Otherwise, let $v_1\in N_1(v_0)$.  By construction, we have that $|N_1(v_1)| \geq |N_1(v_0)|$.  Furthermore, $D$ contains no transitive triangles, so $|N_1(v_0) \cap N_1(v_1)| = 0$.  Thus, $|N_2(v_0)| \geq |N_1(v_1)| \geq |N_1(v_0)|$, and by definition $v_0$ is satisfactory.
\end{proof}

\begin{remark}
Recall that the girth of a undirected graph is the length of its shortest cycle.  Theorem \ref{girth} shows that any counterexample to Conjecture \ref{SeymourConjecture} must have underlying girth of exactly 3.
\end{remark}

\section{Minimal Criminals}

To this point, we have been showing that classes of graphs satisfy Conjecture \ref{SeymourConjecture}.  In this section we reverse course and explore necessary properties of the minimal criminal graphs of $\script{A''}$ from Definition \ref{minimal criminal}.  If Seymour's Second Neighborhood Conjecture is true, then our goal should be to derive such strong constraints on the graphs of $\script{A''}$ that a contradiction is obtained.  On the other hand, if the conjecture is false, then our goal is to find necessary or sufficient conditions for a graph to be in $\script{A''}$; we provide a number of necessary conditions here.

\begin{theorem}
\label{minimal theorem}
If $\script{M}\in \script{A''}$, we have the following:
\begin{enumerate}
\item \label{strong connected} $\script M$ is strongly connected.
\item For each $u\in V(\script{M})$, $A_s(u)\in\{1,2\}$.
\item \label{lots of paths} For every edge $e=(u,v)\in E(\script{M})$, there exists a path of length 1 or 2 avoiding $e$ from $u$ to all but at most 1 element of $\{v\}\cup N_1(v)$.
\item Every edge of $\script{M}$ is the base of either a transitive triangle or a 2-directed diamond.
\item Suppose that $e=(u,v)\in E(\script{M})$ and $|N_1(u)| \leq |N_1(v)|$.  Then $e$ must be the base of at least $|N_1(v)| - |N_1(u)| + 1 $ transitive triangles and the base of at least $|N_1(v)| - |N_1(u)| + 1$ 2-directed diamonds.
\item For any vertex $u\in V(\script{M})$, there exists a vertex $v\in N_{-1}(u)$ such that $A_s(v) = 1$.
\item There exists a directed cycle in $\script{M}$ such that every vertex on the cycle has anti-satisfaction of exactly 1.
\end{enumerate}
\end{theorem}

\begin{proof}

\textbf{Proof of 1:}  Recall that a directed graph is strongly connected if there exists a directed path between any two of its vertices.  Pick an arbitrary vertex $u$ from the vertex set of $\script{M}$.  Now consider $\script{M}' = \script{M}[W(u)]$.  We now pick an arbitrary vertex $v \in W(u)$.  Clearly, $N_{1,\script{M}} (v)\subseteq W(u)$ and $N_{2,\script{M}}(v)\subseteq W(u)$. But this implies that $A_{s,\script{M}'} = |N_{1,\script{M}'} (v)| - |N_{2,\script{M}'} (v)| = |N_{1,\script{M}} (v)| - |N_{2,\script{M}} (v)| = A_{s,\script{M}}$, and hence $v$ is satisfactory in $\script{M}'$ if and only if $v$ is satisfactory in $\script{M}$.  Since by definition $\script{M}$ contains no satisfactory vertices, $v$ cannot be satisfactory in $\script{M}'$.  Thus $\script{M}'$ contains no satisfactory vertices.  But $\script{M}'$ is a subgraph of $\script{M}$, and so by minimality of $\script{M}$ we have that $\script{M} = \script{M}'$.

\textbf{Proof of 2:}  Fix $u$ and pick an arbitrary edge $e = (u,v)\in E(\script{M})$.  Consider the directed graph $Z$ obtained by deleting $e$ from $\script{M}$.  Since $Z$ has fewer edges than $\script{M}$, we have that $Z$ contains a satisfactory vertex.  For each vertex $w\in V(\script{M})$, we note that $|N_{1,Z}(w)| = |N_{1,\script{M}}(w)|$ unless $w=u$, in which case $|N_{1,Z}(u)| = |N_{1,\script{M}}(u)|-1$.  Furthermore, we have that $|N_{2,Z}(w)| \leq |N_{2,\script{M}}(w)|$, except if $w=u$, in which case we have that $|N_{2,Z}(u)| \leq |N_{2,\script{M}}(u)|+1$.  (See Figure \ref{fig:edge deleting}.)

\begin{figure}[tb]
\centering
\includegraphics{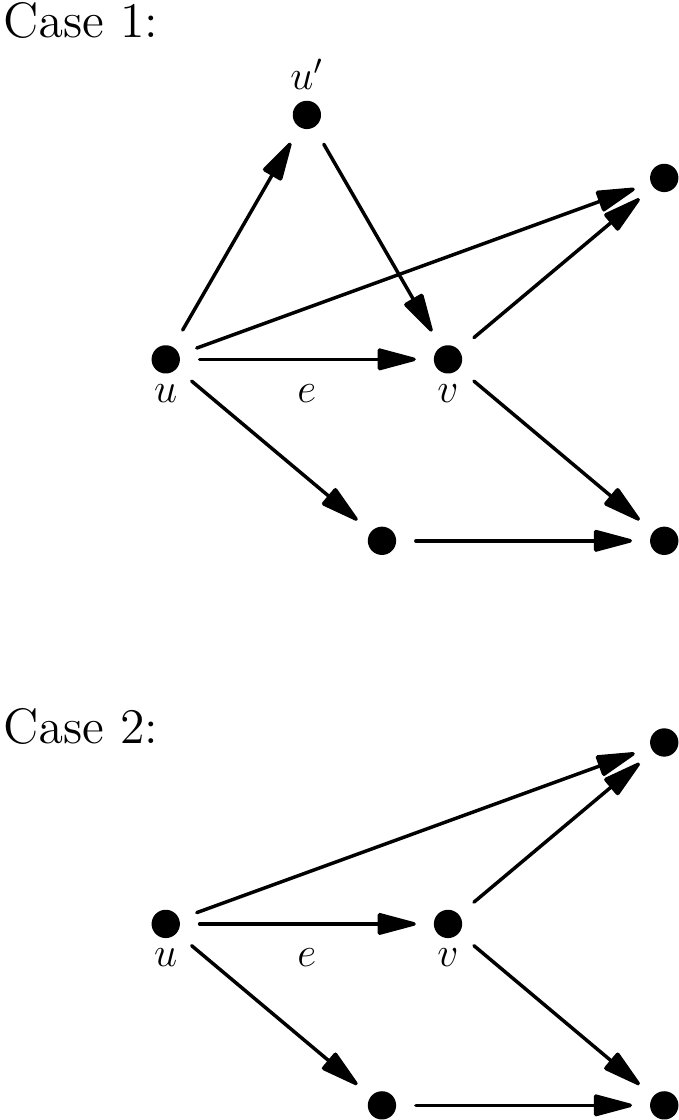}
\caption{Two possible cases resulting from deleting an edge from $\script{M}$.  In Case 1, there is a length 2 path from $u$ to $v$, while in Case 2 no such path exists.  Note that it is possible that deleting $e$ will increase the size of $u$'s second neighborhood, as shown in Case 1.}
\label{fig:edge deleting}
\end{figure}

Thus, we obtain that in $Z$ for $w\neq u\in V(Z)$, $A_{s, Z}(w)\geq A_{s, \script{M}}(w)$, and hence all vertices in $Z$ besides $u$ are not satisfactory.  Thus by process of elimination we have that $u$ is satisfactory in $Z$.  Thus $0 \geq A_{s, Z}(u) = |N_{1,Z}(u)| - |N_{2,Z}(u)| \geq (|N_{1,\script{M}}(u)| - 1) - (|N_{2,\script{M}}(u)|+1)$, and hence we have that $0 < A_{s,\script{M}}(u) = |N_{1,\script{M}}(u)| - |N_{2,\script{M}}(u)| \leq 2$.  Result 2 follows immediately. 

\textbf{Proof of 3:}  We see that $|N_{2,Z}(u)|\geq |N_{2,\script{M}}(u)|$, since otherwise $A_{s, Z}(u)\leq 0$ and $u$ is not satisfactory in $Z$, a contradiction.  Consider now $X = N_{2,Z}(u)\setminus N_{2,\script{M}}(u)$.  We note that $X \subseteq \{v\}$, since $v$ is the only vertex that could have been added to $u$'s second neighborhood in $Z$ (Case 1 in Figure \ref{fig:edge deleting}).  Thus we see that $|N_{2,\script{M}}(u)\setminus N_{2,Z}(u)| \leq 1$, with equality only if $v\in N_{2,Z}(u)$.  

Note that $N_{1,\script{M}}(v) \subseteq N_{1,\script{M}}(u) \cup N_{2,\script{M}}(u)$.  Let $Y = N_{1,\script{M}}(u) \cap N_{1,\script{M}}(v)$ and $Z = N_{2,\script{M}}(u) \cap N_{1,\script{M}}(v)$.  For $y\in Y$, we clearly have a path of length 1 from $u$ to $y$ avoiding $e$ (namely the edge $(u,y)$).  If $|N_{2,\script{M}}(u)\setminus N_{2,Z}(u)| = 0$, then for $z\in Z$, we therefore have a path of length 2 from $u$ to $z$ in $Z$, and considering this path in $\script{M}$ yields a path from $u$ to $z$  avoiding $e$.  And finally, if $|N_{2,\script{M}}(u)\setminus N_{2,Z}(u)| = 1$, then we have a path of length 2 from $u$ to $z$ in $Z$ for all but 1 vertex in $Z$, and as before we have a corresponding path from $u$ to $z$ avoiding $e$.  But in this case, there is a path of length 2 from $u$ to $v$ avoiding $e$, and hence we have obtained the desired result.

\textbf{Proof of 4:} Paths of length 1 from $u$ to $v' \in N_1(v)$ yield transitive triangles with $e$ as the base, and paths of length 2 from $u$ to $v' \in \{v\}\cup N_1(v)$ yield 2-directed diamonds with $e$ as one of the bases.  By part 3, at least one of these structures exists, and hence we are done.

\textbf{Proof of 5:} Since $N_1(v) \setminus (N_1(u) \cap N_1(v)) \subseteq N_2(u)$, we have that $|N_2(u)| \geq |N_1(v)| - |N_1(u) \cap N_1(v)|$.  But since $\script{M}$ contains no satisfactory vertices, we have that $|N_2(u)| < |N_1(u)|$.  By transitivity, we obtain $|N_1(v)| - |N_1(v)\cap N_1(u)| < |N_1(u)|$.  It then follows that $|N_1(v)| - |N_1(u)| < |N_1(v)\cap N_1(u)|$, but $|N_1(v) \cap N_1(u)|$ is the number of transitive triangles having base $e$, so we have proved the first half of part 5. 

To prove the second half of this part, we consider the following cases:

\textit{Case 1}: Suppose there exists a vertex $u'$ such that $(u,u'), (u',v) \in E(\script{M})$.  By part \ref{lots of paths}, we know that $u$ must be connected to at least $|N_1(v)| - 1$ elements of $N_1(v)$ via a path of length 1 or 2 avoiding $e$.  But we see that $u$ is adjacent to at most $|N_1(u) - 2|$ vertices in $N_1(v)$.  Subtracting, we see that $u$ is connected via a path of length 2 avoiding $e$ to at least $|N_1(v)|-1 -(|N_1(u)| -2) = (|N_1(v)| - |N_1(u)|) + 1$ vertices in $N_1(v)$; each of which yields a 2-directed diamond of which $e$ is the base, which is the desired result.

\textit{Case 2}: Suppose there is no such $u'$.  Then again applying part \ref{lots of paths}, it must be that there exists a path of length 1 or 2 avoiding $e$ to each vertex in $N_1(v)$.  But $u$ is adjacent to at most $|N_1(u)| - 1$ of these vertices, and as before we count that there is a path of length 2 avoiding $e$ from $u$ to at least $|N_1(v)| - (|N_1(u)| - 1) = |N_1(v)| - |N_1(u)| + 1$ vertices in $|N_1(v)|$.  Since each of these paths yield a 2-directed diamond with $e$ as the base, we are done.

\textbf{Proof of 6:} In $\script{M}$, pick an arbitrary vertex $u$.  Delete this vertex (and all edges incident with it) and label the resulting directed graph $Z$.  Then in a similar manner to before, one of the vertices in $N_{-1,\script{M}} (u)$ must be satisfactory in $Z$ by vertex minimality of $\script{M}$.  Label this vertex $t$. Since $|N_{1,Z}(t)| = |N_{1,\script{M}}(t)|-1$, $t$ is satisfactory, and $|N_{2,Z}(t)| \subseteq |N_{2,\script{M}}(t)|$ (note that in contrast to deleting an edge, deleting a vertex does not allow any vertices to add vertices to their second neighborhoods), we see that we must have $|N_{2,Z}(t)| = |N_{2,\script{M}}(t)|$.  It is then necessary that $A_{s,\script{M}}(t) = 1$.  Since $u$ was arbitrary, we have obtained the desired result.

\textbf{Proof of 7:}  We apply the same technique as we used Theorem \ref{D.Cycles}.  We present a brief sketch of our proof: by part 5, each vertex in $\script{M}$ has an in-neighbor having anti-satisfaction of exactly 1.  If we begin at an arbitrary vertex and choose one of its in-neighbors having anti-satisfaction of exactly 1, do the same for the resulting vertex, and iterate this process, at some point we must arrive back at a vertex we have already visited, thus constructing a directed cycle of vertices having anti-satisfaction exactly 1.
\end{proof}

Finally, we show that there is not a finite nonzero number of strongly-connected counterexamples to the conjecture.  That is, either the conjecture is true, or there are an infinite number of (non-isomorphic) strongly-connected graphs that violate Conjecture \ref{SeymourConjecture}.  This is especially interesting in light of Part \ref{strong connected} of Theorem \ref{minimal theorem}, which shows that all minimal criminals are strongly connected.

\begin{theorem}
If Seymour's Second Neighborhood Conjecture is false, there are infinitely many non-isomorphic strongly-connected counterexamples to Seymour's Second Neighborhood Conjecture.
\end {theorem}

\begin{proof}
Suppose that Seymour's Second Neighborhood Conjecture is false, and suppose that digraph $D$ is any strongly-connected counterexample to Seymour's Second Neighborhood Conjecture.  (By Part \ref{strong connected} of Theorem \ref{minimal theorem}, such a $D$ must exist.)  Let $H$ be any digraph satisfying the condition $\forall v \in V(H), A_s(v) \geq 0$; that is, all of $H$'s vertices have nonnegative anti-satifaction.  Note that any dicycle satisfies the relevant condition, and hence there exists a choice of $H$ on any number $n$ of vertices, $n\geq 3$.

We now construct a graph $D'$ on $|V(D)|\cdot |V(H)|$ vertices such that $D'$ is a counterexample to Seymour's Second Neighborhood Conjecture, thus proving our theorem.  We define our graph $D'$ as follows:

\begin{itemize}
\item $V(D') = V(D)\times V(H)$
\item If $u = (d_1,h_1), v = (d_2,h_2)\in V(D')$, then $(u,v)\in E(D')$ if and only if either
\begin{enumerate}
\item $d_1=d_2$ and $(h_1, h_2)\in E(H)$, or
\item $d_1 \neq d_2$ and $(d_1, d_2)\in E(D)$.
\end{enumerate}
\end{itemize}

\begin{figure}[tb]
\centering
\includegraphics{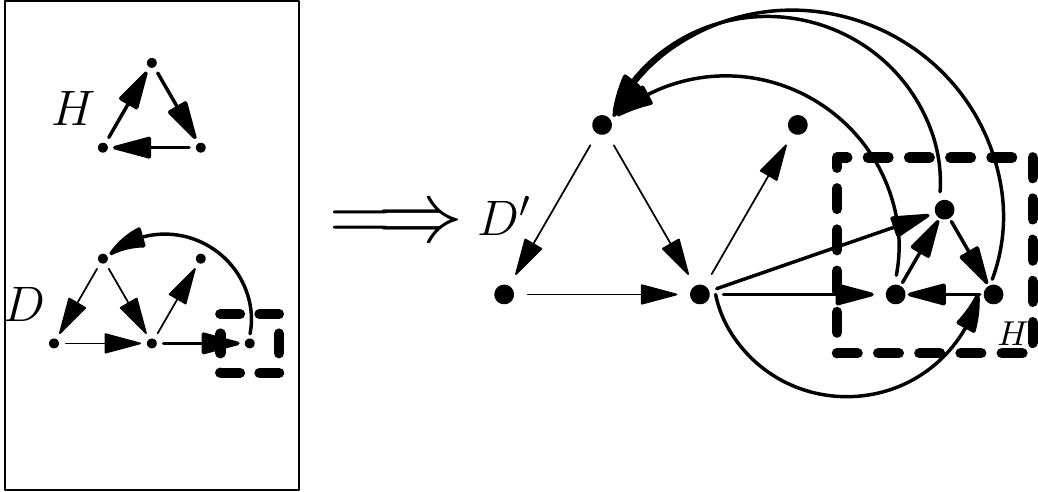}
\caption{A partial representation of the graph $D'$, given $D$ and $H$.  We can think about $D'$ as being made by replacing each vertex of $D$ with a copy of $H$.  Note that for clarity we replace only one vertex in the above picture.}
\label{fig:exploding}
\end{figure}

For any vertex $v = (d, h)\in V(D')$, we calculate that
$$|N_{1, D'} (v)| = |N_{1, H}(h)| + |V(H)| \cdot |N_{1, D'}(d)|,$$
by construction.  Furthermore, we have that
$$|N_{2, D'} (v)| = |N_{2, H}(h)| + |V(H)| \cdot |N_{2, D'}(d)|.$$
We then calculate that
$$\begin{array}{lll}
A_{s,D'} (v) & = & |N_{1, D'} (v)| - |N_{2, D'} (v)|\\
& =  & (|N_{1, H}(h)| - |N_{2, H}(h)|) + |V(H)| (|N_{1, D'}(d)| - |N_{2, D'}(d)|).
\end{array}$$
But by our choice of $H$, we have that $|N_{1, H}(h)| - |N_{2, H}(h)|\geq 0$, and by our choice of $D$ we have that $|N_{1, D'}(d)| - |N_{2, D'}(d)| > 0$.  Hence we obtain $A_{s,D'} (v) > 0$, thus implying that every vertex in $D'$ has positive anti-satisfaction.

Furthermore, $D'$ is strongly connected: fix $(d_1, h_1),(d_2,h_2)\in V(D')$.  If $d_1\neq d_2$, let $d_1, \delta_1, \ldots, \delta_i, d_2$ define a directed path in $D$ from $d_1$ to $d_2$.  Then
$$(d_1, h_1), (\delta_1, h_2),\ldots, (\delta_i, h_2), (d_2, h_2)$$
defines a directed path in $D'$ from $(d_1, h_1)$ to $(d_2,h_2)$.  If $d_1 = d_2$, let $d_3\in N_{1, D}(d_1)$; we know that $(d_1, h_1), (d_3, h_2)$ are adjacent in $D'$, and since $d_2 \neq d_3$ there is a path from $(d_3, h_2)$ to $(d_2, h_2)$ in $D'$, the existence of a path from $(d_1, h_1)$ to $(d_2,h_2)$ follows.  

By definition, we then have that $D'$ is a strongly-connected counterexample to Seymour's Second Neighborhood Conjecture.
\end{proof}

\section{Acknowledgements}

This work was done at the East Tennessee State University REU, NSF grant 0552730, under the supervision of Dr. Anant Godbole.

\bibliographystyle{amsplain}
\bibliography{Seymourbib}

\begin{tabular}{p{3in}l}
\small\textsc{James N. Brantner} & \small\textsc{Greg Brockman}\\[-5pt]
\small\textsc{Erskine College} &\small\textsc{Harvard University}\\[-5pt]
\small\textsc{Due West, SC}&\small\textsc{Cambridge, MA}\\[-5pt]
\small\textsc{United States}&\small\textsc{United States}\\[-5pt]
\small \verb|jbrantne@erskine.edu|&\small\verb|gbrockm@fas.harvard.edu|
\end{tabular}
\vspace{.3in}

\begin{tabular}{p{3in}l}
\small\textsc{Bill Kay} & \small\textsc{Emma Snively}\\[-5pt]
\small\textsc{University of South Carolina} &\small\textsc{Rose-Hulman Institute of Technology}\\[-5pt]
\small\textsc{Columbia, SC}&\small\textsc{Terre Haute, IN}\\[-5pt]
\small\textsc{United States}&\small\textsc{United States}\\[-5pt]
\small \verb|kayw@mailbox.sc.edu|&\small\verb|snivelee@rose-hulman.edu|
\end{tabular}\\
\end{document}